\documentclass{SCIS2018}

\usepackage{colortbl}

\begin{document}

\ArticleType{RESEARCH PAPER}
\Year{2018}
\Month{}
\Vol{61}
\No{}
\DOI{}
\ArtNo{}
\ReceiveDate{}
\ReviseDate{}
\AcceptDate{}
\OnlineDate{}

\title{Optimal Control with \emph{Irregular} Performance}{Optimal Control with \emph{Irregular} Performance}
\author[1]{ZHANG Huanshui}{hszhang@sdu.edu.cn. This version has been accepted by SCIENCE CHINA Information Sciences.}
\author[1]{XU Juanjuan}{}

\AuthorMark{ZHANG H.}

\AuthorCitation{ZHANG H, XU J}


\address[1]{School of Control Science and Engineering, Shandong University, Jinan, Shandong, {\rm 250061}, P.R. China}


\abstract{
In this paper, we solve the long--standing fundamental problem of irregular linear--quadratic (LQ) optimal control, which has received significant attention since the 1960s.
We derive the optimal controllers via the key technique of finding the analytical solutions to two different forward and backward differential equations (FBDEs).
We give a complete solution to the finite--horizon irregular LQ control problem using a new `two-layer optimization' approach. We also obtain the necessary and sufficient condition for the existence of optimal and stabilizing solutions
in the infinite--horizon case in terms of solutions to two Riccati equations
and the stabilization of one specific system.
For the first time, we explore the essential differences between irregular and standard LQ control, making a fundamental contribution to classical LQ control theory.

We show that irregular LQ control is totally different from regular control as the irregular controller must guarantee the terminal state constraint of $P_1(T)x(T)=0$.}
%
%
%
%


\maketitle


%


\section{Introduction}\label{s1}

~~~~Due to its widespread use in modern engineering, linear--quadratic (LQ) optimal control has become one of the most fundamental problems in the field and has received much
attention in recent years~\cite{Ho, Speyer, chenhanfu, chenlizhou, sunliyong, Bell, clementsanderson, Bellman, anderson, lewis, SCIS1, SCIS2, SCIS3, SCIS4, hszhang, zhangIR, zhangOpConsensus, hszhang1}. Since 1950s, researchers have studied classical LQ control problems
where the controller's weight matrix in the cost function
is positive definite~\cite{Kalman, Letov}. Two main approaches have been developed to handle such problems: the maximum principle~\cite{mp} and dynamic
programming~\cite{bellman}. When the controller's weight matrix is positive definite, a unique optimal solution exists, and can be given in terms of the standard Riccati equation, which is clearly regular
in this case.

When the weight matrix is positive semi-definite (and possibly zero), the LQ problem, called {\em singular control}, has remained a significant challenge since the 1960s. Here, there are three main approaches. The first approach is `Transformation in state space' (see \cite{Gurman, Moore, Williems} and references therein for details). The problem with a zero weight matrix $R$ was studied by Ho~\cite{Ho}, who showed that the problem is solvable if the initial values
are of the form $x_2(0)=C_{21}(0)x_1(0)$. Otherwise, an initial impulse control must be applied. It should also be noted that only sufficient conditions have been given for the general case $R\geq 0$. The second approach is the higher--order maximum principle \cite{Gabasov, Krener, Hoehener, Bonnans, XuZhang}. Here, we note that a singular control cannot be found with this approach if the higher derivatives vanish~\cite{Gabasov}.  Furthermore, particular initial values
must usually be chosen to ensure the problem is solvable. The third approach applies a perturbation~\cite{chenhanfu,sunliyong}. In~\cite{sunliyong}, open--loop solutions for irregular LQ problems were obtained by applying a minimize convergence sequence while closed-loop controllers were derived based on the regularity of the Riccati equation.
In summary, although singular control has been studied for more than 50 years, some fundamental problems still remain to be solved, especially irregular LQ control, which is the essence of singular control.

In this paper, 
we consider the solvability of irregular LQ control problems with both finite and infinite horizons.
One key issue in LQ control problems is to find an optimal controller that minimizes the cost function while also stabilizing the original system.
For traditional LQ control problems, there is always a unique positive--definite solution to the standard algebraic Riccati equation
under stabilizability and observability assumptions. However, these conditions are
too strict for infinite--horizon irregular LQ control problems to be solvable. Here, we first derive the necessary and sufficient condition for solving finite--horizon irregular LQ control problems.
Based on a new ``two--layer optimization" approach, we show that the controller entries
of irregular--LQ controller can be derived from two equilibrium conditions in two different layers, unlike in classical regular LQ control where all these entries are obtained from an equilibrium condition in one layer, based on a regular Riccati equation.
Then, we give the necessary and sufficient condition for the existence of optimal and stabilizing solutions to infinite--horizon irregular LQ control problems.
Specially, we derive the optimal controller in terms of two layers, where the first layer is based on the solutions to two Riccati equations and the second one is designed based on the stabilization of one specific system.

The remainder of this paper is organized as follows. Section \ref{s2} introduces the problems studied in this paper. Sections \ref{s3} and \ref{s4} present the solutions to the finite--horizon and infinite--horizon irregular control problems, respectively. Finally, Section \ref{s5} gives some concluding remarks, and the Appendices provide additional proofs.

Throughout this paper, we use the following notation. First, $R^n$
denotes the family of $n$-dimensional vectors, $x'$ represents the
transpose of $x$, and $\|x\|^2=x'x$. Stating that a symmetric matrix $M>0\ (\geq 0)$ means
it is strictly positive definite (positive semi--definite), and $Range(M)$ represents its range. Finally, $M^{\dag}$ is the Moore--Penrose inverse~\cite{pinv} of the matrix $M$ if it satisfies $MM^{\dag}M=M,~M^{\dag}MM^{\dag}=M^{\dag},
(MM^{\dag})¡¯=MM^{\dag}$, and $(M^{\dag}M)'=M^{\dag}M$.


\section{Problem Formulation}\label{s2}

Consider the linear system governed by the following differential equation:
\begin{eqnarray}
\dot{x}(t)&=&Ax(t)+Bu(t),~x(t_0)=x_0,\label{d1}
\end{eqnarray}
where $x\in R^n$ is the system's state, $u\in R^m$ is its control input, and $x_0$ is its initial value.
The matrices $A$ and $B$ are constant and have compatible dimensions.
The cost function is given by
\begin{eqnarray}
J_T(t_0,x_0;u)=\int_{t_0}^T[x'(t)Qx(t)+u'(t)Ru(t)]dt+x'(T)Hx(T),\label{df}
\end{eqnarray}
where $Q\geq 0$, $R\geq 0$, and $H\geq0$ are symmetric matrices with compatible dimensions.
Here, we consider two problems, namely finite--horizon and infinite--horizon optimal control. The finite--horizon optimal control problem can be stated as follows.

\textbf{Problem 1:} For any given initial pair $(t_0,x_0)$, find a $u^*(t)$ such that
\begin{eqnarray}
J^*_T(t_0,x_0;u^*)=\min_{u(\cdot)}J_T(t_0,x_0;u).\nonumber
\end{eqnarray}

\begin{remark}
Unlike with standard LQ control, the weight matrix $R$ in the cost function (\ref{df}) is singular, which makes the problem significantly more complex, as mentioned in Section \ref{s1}. We will show that the difficulty of singular control problems is caused by the so--called irregularity, that is, the fact that the controller cannot be derived from the equilibrium condition. In other words, irregular control problems cannot be solved by simply completing the square. To make this clearer, we present the following example.

Consider the system
\begin{eqnarray}
\dot{x}(t)&=&\left[
                            \begin{array}{cc}
                              1 & 1 \\
                            \end{array}
                          \right]u(t),~x(t_0)=x_0,\nonumber
\end{eqnarray} and the cost function \begin{eqnarray}
J_T(t_0,x_0;u)=\int_{t_0}^Tu'(t)\left[
                                   \begin{array}{cc}
                                     1 & 0 \\
                                     0 & 0 \\
                                   \end{array}
                                 \right]u(t)dt+x'(T)x(T).\label{de1}
\end{eqnarray}
The solution to the Riccati equation $0=\dot{P}(t)-P^2(t)$ with $P(T)=1$ is $P(t)=-\frac{1}{t-T-1}$.
We then find that $Range \Big(B'P(t)\Big)\nsubseteq Range (R)$, implying that we cannot obtain $u(t)$ from the equation $Ru(t)+B'P(t)x(t)=0$ for any $x(t).$

In addition, taking the derivative of $x'(t)P(t)x(t)$ yields
\begin{eqnarray}
\frac{d}{dt}[x'(t)P(t)x(t)]=2u'(t)B'P(t)x(t)+x'(t)P^2(t)x(t).\nonumber
\end{eqnarray}
Then, by integrating this from $t_0$ to $T$, the cost function (\ref{de1}) can be rewritten as
\begin{eqnarray}
J_T(t_0,x_0;u)&=&x'(t_0)P(t_0)x(t_0)+\int_{t_0}^T\Big(u'(t)Ru(t)+2u'(t)B'P(t)x(t)+x'(t)P^2(t)x(t)\Big)dt\nonumber\\
&=&x'(t_0)P(t_0)x(t_0)+\int_{t_0}^T\Big[\Big(u(t)+R^{\dag}B'P(t)x(t)\Big)'R\Big(u(t)+R^{\dag}B'P(t)x(t)\Big)Ru(t)\nonumber\\
&&+2u'(t)(I-RR^{\dag})B'P(t)x(t)\Big]dt.\nonumber
\end{eqnarray}
From this cost function and the fact that $Range \Big(B'P(t)\Big)\nsubseteq Range (R)$, we can see that the optimal controller cannot be obtained by completing the square.

Given this, we cannot apply any of the classical approaches to this irregular optimal control problem. However, it is nonetheless solvable; in fact, one optimal solution is given by $u(t)=-\left[
               \begin{array}{c}
                 0 \\
                 1 \\
               \end{array}
             \right]\frac{x(t_0)}{T-t_0}$ or $u(t)=\left[
               \begin{array}{c}
                 0 \\
                 1 \\
               \end{array}
             \right]\frac{x(t)}{t-T}$, and the corresponding optimal cost is $0$.

In literature, there exists some work on LQ control in the $R\geq0$ case, even for stochastic control with indefinite $R$ (see~\cite{zhou2} and references therein). However, the problems they studied were assumed to be regular, and the solvability of the irregular LQ problem remains a challenge. One contribution of this paper is that we present the complete solution to the irregular LQ problem.

\end{remark}

\begin{remark}
Irregular LQ control has broad applications in numerous areas, such as engineering systems, economic models, and the production of natural resources (\cite{Bell}, \cite{Kliger, Hsia}). For example, a fundamental space navigation problem is to control rockets so as to minimize the propellant expenditure based on a mathematical model that includes bounds on the rocket's thrust, which is a typical irregular LQ problem~\cite{Bell}.

\end{remark}
In the infinite--horizon case, we are only interested in stabilizing controllers; thus, we introduce the admissible control set
\begin{eqnarray}
\mathcal{U}&=&\Big\{u(t),t\in[0,\infty]\Big|\int_{0}^\infty\|u(t)\|^2dt<\infty, u(t) ~\mbox{is~stabilizing~control}\Big\}.\nonumber
\end{eqnarray}
The infinite--horizon problem is described in detail below.

\textbf{Problem 2:} Find a control $u(t)\in \mathcal{U}$ which minimizes the cost function $J(x_0;u)$ and stabilizes the system (\ref{d1}) where $J(x_0;u)$ is defined by
\begin{eqnarray}
J(x_0;u)=\int_{0}^\infty[x'(t)Qx(t)+u'(t)Ru(t)]dt.\label{d2}
\end{eqnarray}

We now discuss this irregularity in detail. It is well known~\cite{hszhang1} that Problem 1 is solvable if and only if the following forward and backward differential equations (FBDEs) are solvable:
\begin{eqnarray}
\dot{x}(t)&=&Ax(t)+Bu(t),~x(t_0)=x_0,\label{1}\\
\dot{p}(t)&=&-[A'p(t)+Qx(t)],~p(T)=Hx(T),\label{2}\\
0&=&Ru(t)+B'p(t).\label{3}
\end{eqnarray}
This implies that the key to deriving optimal solutions is the solvability of the FBDEs (\ref{1})--(\ref{3}). To this end, we introduce the following Riccati equation:
\begin{eqnarray}
0=\dot{P}(t)+A'P(t)+P(t)A+Q-P(t)BR^{\dag}B'P(t), \label{dr1}
\end{eqnarray}
where the terminal condition is given by $P(T)=H$ and $R^{\dag}$ represents the Moore--Penrose inverse of $R$. Throughout this paper, we assume that
the Riccati equation (\ref{dr1}) is solvable. The solvability condition has been given in Corollary 2.10 of Chapter 6 in \cite{yong}.

In the $R>0$ case, if we replace $p(t)$ with $P(t)x(t)$, the equilibrium condition (\ref{3}) becomes
\begin{eqnarray}
0&=&Ru(t)+B'P(t)x(t).\label{f1}
\end{eqnarray}
Thus, when \begin{eqnarray}
Range \Big(B'P(t)\Big)\subseteq Range \Big(R\Big),\label{d32}
\end{eqnarray}
the linear equation (\ref{f1}) is solvable for any $x(t)\in R^n$, which implies the optimization problem is solvable as well. Conversely, when
\begin{eqnarray}
Range \Big(B'P(t)\Big)\not\subseteq Range \Big(R\Big),\label{d33}
\end{eqnarray}
the linear equation (\ref{f1}) is unsolvable. This case of (\ref{d33}) is termed as irregular control.
As shown in the previous example (Remark 1), an optimal control cannot be obtained from the equilibrium condition or by
reformulating the cost function by completing the square. However, the irregular case may still be solvable. In the next section, we give
a complete solution to the irregular LQ control problem.

In summary, as discussed above, the irregularity of the problem implies that we cannot simply complete the square for the cost function (\ref{df}).

\section{Solution to Problem 1}\label{s3}

In this section, we study the finite--horizon optimal control problem, and obtain results for both the regular and irregular cases.
First, we give the optimal solution in the regular case.
\begin{theorem}\label{lemre}
In the regular case (\ref{d32}), the optimal solution to Problem 1 is given by
\begin{eqnarray}
u(t)&=&-R^{\dag}B'P(t)x(t)+(I-R^{\dag}R)z(t),
\end{eqnarray}
where $z(t)$ is an arbitrary vector of compatible dimension. In this case, the optimal cost is
\begin{eqnarray}
J^*_T(t_0,x_0;u^*)=x_0'P(t_0)x_0.\label{dinf7}
\end{eqnarray}
\end{theorem}
{\em Proof}. This result can be obtained directly from~\cite{zhou2}.

Now, we turn to the irregular case, namely $Range \Big(B'P(t)\Big)\not\subseteq Range \Big(R\Big)$. First, we transform the task of solving Problem 1 into that of solving new FBDEs.
Before presenting the result, we introduce some notations for future convenience. Without loss of generality, let $rank (R)=m_0< m$, so $rank(I-R^{\dag}R)=m-m_0>0$. There is also an elementary row transformation matrix $T_0$ such that
\begin{eqnarray}
T_0(I-R^{\dag}R)=\left[
                              \begin{array}{c}
                                 0  \\
                                 \Upsilon_{T_0}\\
                              \end{array}
                            \right], \label{jnYC1}
\end{eqnarray}
where $\Upsilon_{T_0}\in R^{[m-m_0]\times m}$ is of full row rank. In addition, let
\begin{eqnarray}
A_0(t)&=&A-BR^{\dag}B'P(t),~D_0=-BR^{\dag}B',~\left[ \begin{array}{cc}\ast & G_0\\
                              \end{array}
                            \right] ={T_0}^{-1},\nonumber\\
 \left[ \begin{array}{cc} \ast & C_0'(t) \\
                              \end{array} \right] &=&P(t)B\Big(I-R^{\dag}R\Big){T_0}^{-1},~~
\left[ \begin{array}{cc}\ast & B_0\\
                              \end{array}
                            \right] =B\Big(I-R^{\dag}R\Big){T_0}^{-1},\nonumber
\end{eqnarray}
and define
\begin{eqnarray}
0&=&\dot{P}_1(t)+P_1(t)A_0(t)+A_0'(t)P_1(t)+P_1(t)D_0P_1(t), \label{d3}
\end{eqnarray}
where the terminal value $ P_1(T)$ is yet to be determined. Next, Lemma \ref{lemma-1} states that solving Problem 1 is equivalent to solving new FBDEs.
\begin{lemma}\label{lemma-1}
Given condition (\ref{d33}), Problem 1 is solvable if and only if there exists a $u_1(t)\in R^{m-m_0}$ such that
\begin{eqnarray}
0=C_0(t)x(t)+B_0'\Theta(t),\Theta(T)=0,\label{n2}
\end{eqnarray}  where $u_1(t)$, $x(t)$, and $\Theta(t)$ satisfy the FBDEs
\begin{eqnarray}
\dot{x}(t)&=&A_0(t)x(t)+D_0\Theta(t)+B_{0}u_1(t),\label{n3}\\
\dot{\Theta}(t)&=&-\Big(A'_0(t)\Theta(t)+C'_0(t)u_1(t)\Big).\label{c8}
\end{eqnarray}

\end{lemma}
\emph{Proof.} The proof is given in \ref{apendix1}.

We are now in a position to give the solvability condition for Problem 1.
\begin{theorem}\label{lemir}
Given condition (\ref{d33}), Problem 1 is solvable if and only if there exists a $P_1(t)$ in (\ref{d3}) with terminal value $P_1(T)$ such that
\begin{eqnarray}
0=C_0(t)+B_0'P_1(t), ~~~t_0\leq t\leq T,\label{d4}
\end{eqnarray}
and a $u_{1}(t)$ that achieves
\begin{eqnarray}
P_1(T)x(T)=0,\label{d5}
\end{eqnarray} where $x(t)$ obeys
\begin{eqnarray}
\dot{x}(t)&=&\Big(A_0(t)+D_0P_1(t)\Big)x(t)+B_0u_1(t)\label{d24}
\end{eqnarray}
and has initial value $x(t_0)=x_0$. In this case, the optimal controller $u(t)$ is given by
\begin{eqnarray}
u(t)=-R^{\dag}B'\Big(P(t)+P_1(t)\Big)x(t)+G_0u_1(t),\label{df1}
\end{eqnarray}
and the optimal cost is given by
\begin{eqnarray}
J_T^*(t_0,x_0;u^*)=x_0'\Big(P(t_0)+P_1(t_0)\Big)x_0.\nonumber
\end{eqnarray}

\end{theorem}
\emph{Proof.} The proof is given in \ref{apendix2}.

\begin{corollary}
The following statements hold.
\begin{enumerate}
  \item Under the conditions of Theorem~\ref{lemre}, $P(t)\geq 0$.
  \item Under the conditions of Theorem~\ref{lemir}, $P(t)+P_1(t)\geq0$.
\end{enumerate}
\end{corollary}
\emph{Proof.} Since the matrices $Q$, $R$, and $H$ are positive semi--definite, the cost function (\ref{df}) is non--negative. Thus, the optimal costs obtained in Theorems~\ref{lemre} and \ref{lemir} are non--negative.
Considered together with the arbitrariness of the initial time and value, these results follow, completing the proof.

Next, we formulate the open--loop and closed--loop solutions.
\begin{theorem}\label{Th-1}
Given condition (\ref{d33}), Problem 1 is open--loop solvable if and only if there exists a $P_1(t)$ in (\ref{d3}) such that (\ref{d4}) holds and
\begin{eqnarray}
Range \big[P_1(t_0)\big]\subseteq
Range \Big(G_1[t_0,T]\Big),\label{df4}
\end{eqnarray}
where the Gramian matrix $G_1[t_0,T]$ is defined by
\begin{eqnarray}
G_1[t_0,T]=\int_{t_0}^TP_2(t_0,s)C_0'(s)C_0(s)P_2'(t_0,s)ds,
\end{eqnarray}
and $P_2(t,s)$ satisfies
$\dot{P}_2(t,s)=-A_0'(t)P_2(t,s)$, with $P_2(t,t)=I$.
In this case, the open--loop solution is given by (\ref{df1}), while $u_1(t)$ is given by
\begin{eqnarray}
u_1(t)&=&C_0(t)P_2'(t_0,t)G_1^{\dag}[t_0,T]P_1(t_0)x_0.\label{df5}
\end{eqnarray}
\end{theorem}
\emph{Proof.} The proof is given in \ref{apendix3}.

Next, we consider the closed--loop solution. To this end, let
$P_1(t)=\left[
         \begin{array}{cc}
           P^1_{11}(t) & P^1_{12}(t) \\
           {P^1_{12}}'(t) & P^1_{22}(t) \\
         \end{array}
       \right]$.
Assume that $P_1(t)$ is singular, then there exists a $\mathcal{T}_1(t)$ such that
a $\mathcal{T}_1'(t)P_1(t)\mathcal{T}_1(t)=\left[
         \begin{array}{cc}
           \hat{P}(t) & 0 \\
           0 & 0 \\
         \end{array}
       \right]$,
where $\hat{P}(t)$ is invertible. Let $
\dot{\mathcal{T}}_1'(t)\mathcal{T}_1(t)=\left[
      \begin{array}{c}
        \tilde{T}_{1}(t) \\
        \tilde{T}_{2}(t) \\
      \end{array}
    \right], \mathcal{T}_1'(t)B_0\mathcal{T}_1(t)=\left[
      \begin{array}{c}
        B_{1}(t) \\
        B_{2}(t) \\
      \end{array}
    \right]$,
and
$\mathcal{T}_1'(t)\Big(A_0(t)+D_0P_1(t)\Big)\mathcal{T}_1(t)=\left[
                                 \begin{array}{c}
                                   \hat{A}_{1}(t) \\
                                   \hat{A}_{2}(t) \\
                                 \end{array}
                               \right]$.

\begin{theorem}\label{thd3}
Given condition (\ref{d33}), assuming that (\ref{d4}) holds and that there exists a $K(t)$ such that
\begin{eqnarray}
\tilde{T}_{1}(t)+\hat{A}_{1}(t)+B_{1}(t)\mathcal{T}_1'(t)K(t)\mathcal{T}_1(t)=\left[
     \begin{array}{cc}
       \frac{I}{t-T} & 0 \\
     \end{array}
   \right],\label{df7}
\end{eqnarray}
then Problem 1 is closed--loop solvable.
The closed--loop solution is given by (\ref{df1}) with
$u_1(t)=K(t)x(t)$, where $K(t)$ satisfies (\ref{df7}).
\end{theorem}
\emph{Proof.} The proof is given in \ref{apendix4}.

\section{Solution to Problem 2}\label{s4}

In this section, we focus on the infinite--horizon optimal control problem, giving necessary and sufficient conditions for the solvability of this optimization problem in both the regular and irregular cases.
In the discussion below, we make use of the following algebraic Riccati equation:
\begin{eqnarray}
0=A'P+PA+Q-PBR^{\dag}B'P. \label{dinf4}
\end{eqnarray}

\subsection{Regular Case}

First, we consider the case where $Range \Big(B'P(t)\Big)\subseteq Range \Big(R\Big)$.

\begin{theorem}\label{th1}
Given the condition (\ref{d32}), there exists an optimal and stabilizing solution to Problem 2
if and only if the following conditions hold.
\begin{enumerate}
  \item There exists a solution $P$ to (\ref{dinf4}) such that $P\geq 0$.
  \item The system $(A-BR^{\dag}B'P,B(I-R^{\dag}R))$ is stabilizable.
\end{enumerate}
In this case, the optimal and stabilizing solution is given by
\begin{eqnarray}
u(t)&=&-R^{\dag}B'Px(t)+(I-R^{\dag}R)z(t),\label{r15}\\
z(t)&=&Kx(t),\label{r16}
\end{eqnarray}
where $K$ must be chosen such that the matrix
\begin{eqnarray}
A-BR^{\dag}B'P+B(I-R^{\dag}R)K\label{r19}
\end{eqnarray}
is stable.
\end{theorem}
\emph{Proof.} The proof can be found in~\cite{zhangOpConsensus}.

\begin{remark}
Without loss of generality, we can assume the existence of an orthogonal matrix $T$ such that $T'RT=\left[
                                   \begin{array}{cc}
                                     R_1 & 0 \\
                                     0 & 0 \\
                                   \end{array}
                                 \right]
$, with $R_1>0$. If we let $\bar{u}(t)=T'u(t)$ and $BT=\left[
                                           \begin{array}{cc}
                                             B_1 & B_2 \\
                                           \end{array}
                                         \right]
$, the system (\ref{d1}) reduces to
\begin{eqnarray}
\dot{x}(t)=Ax(t)+\left[
                                           \begin{array}{cc}
                                             B_1 & B_2 \\
                                           \end{array}
                                         \right]\bar{u}(t).\nonumber
\end{eqnarray}
The Riccati equation (\ref{dinf4}) then becomes
\begin{eqnarray}
0=A'P+PA+Q-PB_1R_1^{-1}B_1'P. \label{d34}
\end{eqnarray}
By applying the standard results in \cite{lewis} and \cite{hszhang}, we have that when $(A,Q)$ is observable and $(A ,B_1)$ is stabilizable, the algebraic
Riccati equation (\ref{d34}) admits a unique positive--definite solution $P$. An iterative algorithm for solving (\ref{d34}) is given in Corollary 1 of \cite{ZhangFangfang}.
In addition, the matrix $A-B_1R_1^{-1}B_1'P$ is stable. In this case, $K$ in (\ref{r16}) can be chosen to be the zero matrix, leading to the matrix in (\ref{r19}) being exactly $A-BR^{\dag}B'P=A-B_1R_1^{-1}B_1'P$, which is stable.
Thus, if $(A,Q)$ is observable and $(A,B_1)$ is stabilizable, then there exists an optimal and stabilizing solution that
minimizes the cost function (\ref{d2}), according to Theorem~\ref{th1}. Accordingly, these conditions are sufficient for such a solution to exist.
However, the reverse is not true. For example, consider the system
$\dot{x}(t)=\left[
               \begin{array}{cc}
                 1 & 1 \\
               \end{array}
             \right]\left[
                      \begin{array}{c}
                        u_1(t) \\
                        u_2(t) \\
                      \end{array}
                    \right]$
and the cost function $
J=\int_0^\infty \left[
               \begin{array}{cc}
                 u_1'(t) & u_2'(t) \\
               \end{array}
             \right]\left[
                      \begin{array}{cc}
                        1 & 0 \\
                        0 & 0 \\
                      \end{array}
                    \right]\left[
                      \begin{array}{c}
                        u_1(t) \\
                        u_2(t) \\
                      \end{array}
                    \right]dt$.
According to Theorem~\ref{th1}, the optimal controller is
$u(t)=\left[
                      \begin{array}{c}
                        0 \\
                        u_2(t) \\
                      \end{array}
                    \right]$,
where $u_2(t)$ can be chosen as
$u_2(t)=-x(t)$,
stabilizing the state $x(t)$, and the optimal cost function is
$J^*=0$.
However, in this case we find that
$A=0$, $Q=0$, and $B_1=1,$
indicating that $(A,Q)$ is unobservable. Thus, it is not necessary for $(A,Q)$ to be observable and $(A ,B_1)$ to be stabilizable.
\end{remark}

\subsection{Irregular Case}

Next, we consider the case where $Range \Big(B'P(t)\Big)\not\subseteq Range \Big(R\Big).$ Here, we define the algebraic Riccati equation corresponding to (\ref{d3}) as
\begin{eqnarray}
0&=&P_1A_0+A_0'P_1+P_1D_0P_1, \label{d6}
\end{eqnarray}
where $A_0=A-BR^{\dag}B'P$. In addition, let $\left[ \begin{array}{cc} \ast & C_0' \\
                              \end{array} \right] =PB(I-R^{\dag}R){T_0}^{-1}$.

\begin{theorem}\label{infir}
Given condition (\ref{d33}), there exists an optimal and stabilizing solution to Problem 2 if and only if the following conditions hold.
\begin{enumerate}
  \item There exist solutions $P$ and $P_1$ to (\ref{dinf4}) and (\ref{d6}), respectively, that satisfy $P+P_1\geq 0$.
  \item $P$ and $P_1$ satisfy
  \begin{eqnarray}
  C_0+B_0'P_1=0.\label{d10}
  \end{eqnarray}
  \item The system $(A_0+D_0P_1,B_0)$ is stabilizable, i.e., there exists a controller $u_1(t)$ such that the following system is stable:
  \begin{eqnarray}
\dot{x}(t)&=&\Big(A_0+D_0P_1\Big)x(t)+B_0u_1(t).\label{d8}
\end{eqnarray}
\end{enumerate}
In this case, the optimal controller is given by
\begin{eqnarray}
u(t)=[-R^{\dag}B'(P+P_1)+G_0K]x(t),\label{df8}
\end{eqnarray}
where $K$ must be chosen such that $A_0+D_0P_1+B_0K$ is stable,
and the optimal cost is
\begin{eqnarray}
J^*=x_0'(P+P_1)x_0.\nonumber
\end{eqnarray}
\end{theorem}
\emph{Proof.}
The proof is given in \ref{apendix5}.

\section{Conclusions}\label{s5}

In this paper, we have solved the long--standing fundamental problem of optimal LQ control with irregular performance. We have obtained optimal solutions to this problem for both the
finite and infinite--horizon cases. In particular, for the finite--horizon case, we have been able to obtain a complete solution by applying a `two-layer optimization' approach.
For the infinite--horizon case, we have given a necessary and sufficient
condition for the existence of an optimal and stabilizing solution in terms of solutions to two Riccati equations
and the stabilization of one specific system. These results have possible applications in areas such as irregular measurement feedback control, robust control, $H_{\infty}$ control, and stochastic control.

We would also like to emphasize that the proposed technique used in this paper, namely the solving of an associated FBDEs, is a very general strategy, with which we have successfully solved other complicated control problems such as stochastic LQ control with time delay~\cite{hszhang,hszhang1}, Stackerberge game control~\cite{game}, mean field stochastic control~\cite{meanfield},
stabilization of NCSs with simultaneous transmission delay and packet dropout~\cite{ncs}.


\Acknowledgements{This work was supported by the National Natural Science Foundation of China (Grant Nos. 61633014, 61573221, and 61873332) and
the Qilu Youth Scholar Discipline Construction Funding from Shandong University.}



\begin{appendix}

\section{Proof of Lemma~\ref{lemma-1}}\label{apendix1}

\emph{Proof of Necessity.} Based on the discussion of (\ref{f1})-(\ref{d33}), we can see that $p(t)\neq P(t)x(t)$ under the condition (\ref{d33}),
where $P(t)$ is the solution to (\ref{dr1}). We therefore define a new variable $\Theta(t)$ as
\begin{eqnarray}
p(t)=P(t)x(t)+\Theta(t),\label{21}
\end{eqnarray}
where it is clear that $\Theta(T)=0$.
Next, we aim to derive the new FBDEs (\ref{n2})-(\ref{c8}) under the solvability of Problem 1.

First, we take the derivative of (\ref{21}), obtaining
\begin{eqnarray}
\dot{p}(t)=\dot{P}(t)x(t)+P(t)\big[Ax(t)
+Bu(t)\big]+\dot{\Theta}(t).\label{zz4}
\end{eqnarray}
From (\ref{2}) and (\ref{21}), we then find that
\begin{eqnarray}
\dot{p}(t)=-\big[A'P(t)x(t)+A'\Theta(t)
+Qx(t)\big].\label{11}
\end{eqnarray}
By comparing (\ref{zz4}) and (\ref{11}), we obtain
\begin{eqnarray}
0&=&\dot{P}(t)x(t)+P(t)Ax(t)+P(t)Bu(t)+\dot{\Theta}(t)
+A'P(t)x(t)+A'\Theta(t)+Qx(t).\label{r6}
\end{eqnarray}
Second, we aim to find the controller $u(t)$ and the new equilibrium condition (\ref{n2}).
By using (\ref{21}), we can formulate the equilibrium condition (\ref{3}) as
\begin{eqnarray}
0=Ru(t)+B'p(t)=Ru(t)+B'P(t)x(t)+B'\Theta(t).\label{25}
\end{eqnarray}
Taken together with (\ref{d33}), this can also be written as
\begin{eqnarray}
u(t)&=&-R^{\dag}\Big(B'P(t)x(t)+B'\Theta(t)\Big)+\Big(I-R^{\dag}R\Big)z(t),\label{n1}
\end{eqnarray}
where $z(t)$ is an arbitrary vector with compatible dimension such that the following equality holds:
\begin{eqnarray}
0&=&\Big(I-RR^{\dag}\Big)\Big(B'P(t)x(t)+B'\Theta(t)\Big). \label{nz1}
\end{eqnarray}
Let \begin{eqnarray}
T_0\Big(I-R^{\dag}R\Big)z(t)=\left[
                              \begin{array}{c}
                                 0  \\
                                 u_1(t)\\
                              \end{array}
                            \right], \label{jn1}
\end{eqnarray}
where $u_1(t)=\Upsilon_{T_0} z(t)\in R^{m-m_0(t)}$.
Now, we can rewrite (\ref{nz1}) as (\ref{n2}). Note that
\begin{eqnarray}
I-RR^{\dag}
&=&(I-RR^{\dag})(I-RR^{\dag})\nonumber \\
&=&(I-RR^{\dag})T'_0(T^{-1}_0)'(I-RR^{\dag})\nonumber \\
&=&\left[
\begin{array}{cc}
                                 0 &
                                 \Upsilon'_{T_0}\\
                              \end{array}
                            \right](T^{-1}_0)'(I-RR^{\dag}),
\end{eqnarray}
where we have used (\ref{jnYC1}) to derive the last equality. By using the definitions below equation (\ref{jnYC1}), we can rewrite (\ref{nz1}) as
\begin{eqnarray}
0&=&\Upsilon'_{T_0}\Big[C_0(t)x(t)+B_0'\Theta(t)\Big].\label{nYc2}
\end{eqnarray}
Note that $\Upsilon'_{T_0}$ is of full column rank, and thus (\ref{nYc2}) can be directly rewritten as (\ref{n2}).

Third, we derive the dynamic of $\Theta(t)$. Substituting (\ref{n1}) into (\ref{r6}) and using (\ref{dr1}) yields
\begin{eqnarray}
0&=&\dot{P}(t)x(t)+P(t)Ax(t)
+A'P(t)x(t)+A'\Theta(t)+Qx(t)+\dot{\Theta}(t)\nonumber\\
&&-P(t)BR^{\dag}\Big(B'P(t)x(t)+B'\Theta(t)\Big)+P(t)B(I-R^{\dag}R)z(t)\nonumber\\
&=&\dot{\Theta}(t)+\Big(A'-P(t)BR^{\dag}B'\Big)\Theta(t)+P(t)B(I-R^{\dag}R)z(t).\label{ZJ1}
\end{eqnarray}
Since $(I-R^{\dag}R)^2=I-R^{\dag}R$, we find that
\begin{eqnarray}
P(t)B(I-R^{\dag}R)z(t)
&=&P(t)B(I-R^{\dag}R)T_0^{-1}T_0(I-R^{\dag}R)z(t)\nonumber\\
&=&P(t)B(I-R^{\dag}R)T_0^{-1}\left[
                              \begin{array}{c}
                                 0  \\
                                 u_1(t)\\
                              \end{array}
                            \right]\nonumber\\
&=&\left[
     \begin{array}{cc}
       * & C_0'(t) \\
     \end{array}
   \right]\left[
                              \begin{array}{c}
                                 0  \\
                                 u_1(t)\\
                              \end{array}
                            \right]=C_0'(t)u_1(t).\label{YCJ1}
\end{eqnarray}
Thus, from (\ref{ZJ1}) we obtain
$\dot{\Theta}(t)=-\Big[A'_0(t)\Theta(t)+C'_0(t) u_1(t)\Big],$
which implies that the dynamics of $\Theta(t)$ are given by (\ref{c8}).

Finally, we derive the dynamics equation (\ref{n3}). By substituting (\ref{n1}) into (\ref{1}) and combining this with the fact that
$B(I-R^{\dag}R)z(t)=B_0u_1(t)$, which can be obtained in a similar way to (\ref{YCJ1}),
we can derive the state dynamics (\ref{n3}).

\emph{Proof of Sufficiency.} Now, we show Problem 1 is solvable if there exists a $u_1(t)$ that enables us to achieve (\ref{n2}). In fact, if (\ref{n2}) is true then (\ref{n1}) and (\ref{nz1}) can be jointly rewritten as (\ref{25}).
Further, by reversing the process for (\ref{21})--(\ref{25}), we can easily verify that
$p(t)=P(t)x(t)+\Theta(t)$, where $x(t)$ and $\Theta(t)$
satisfy (\ref{n2})--(\ref{c8}), solves (\ref{1})--(\ref{3}). Thus, Problem 1 is solvable, completing the proof.

\section{Proof of Theorem~\ref{lemir}}\label{apendix2}

\emph{Proof of Sufficiency.} Based on Lemma~\ref{lemma-1}, it is sufficient to verify that $(\Theta(t),x(t))=(P_1(t)x(t),x(t))$ is the solution to the FBDEs (\ref{n2})--(\ref{c8}).
Taking the derivative of $P_1(t)x(t)$ yields
\begin{eqnarray}
\frac{d[P_1(t)x(t)]}{dt}
&=&\dot{P}_1(t)x(t)+P_1(t)[A_0(t)+D_0P_1(t)] x(t)+P_1(t)B_0u_1(t)\nonumber\\
&=&-A_0'(t)P_1(t)x(t)+P_1(t)B_0u_1(t)\nonumber\\
&=&-A_0'(t)P_1(t)x(t)-C_0'(t)u_1(t),\label{df2}
\end{eqnarray}
where we have used (\ref{d3}) and (\ref{d4}) to derive the last equality.
In addition, again using (\ref{d4}), we have
\begin{eqnarray}
C_0(t)x(t)+B_0'P_1(t)x(t)=0.\label{df3}
\end{eqnarray}
By comparing (\ref{n2}), (\ref{n3}), and (\ref{c8}) with (\ref{df3}), (\ref{df2}), and (\ref{d24}), we can see that (\ref{n2})--(\ref{c8})
is solvable  with
$\Theta(t)=P_1(t)x(t)$ if $P_1(T)x(T)=0$. Thus, based on Lemma~\ref{lemma-1}, Problem 1 is solvable.

\emph{Proof of Necessity.} This proof is divided into two parts. First, we consider the case where the optimal solution is of closed--loop form, namely $u_1(t)=K_1(t)x(t)$. Based on Lemma~\ref{lemma-1},
(\ref{n2})--(\ref{c8}) are solvable if Problem 1 is solvable. By substituting $u_1(t)=K_1(t)x(t)$ into (\ref{n3})-(\ref{c8}), we obtain
\begin{eqnarray}
\dot{x}(t)&=&A_0(t)x(t)+D_0\Theta(t)+B_{0}K_1(t)x(t),\nonumber\\
\dot{\Theta}(t)&=&-\Big(A'_0(t)\Theta(t)+C'_0(t)K_1(t)x(t)\Big).\nonumber
\end{eqnarray}
Solving the above FBDEs gives us $\Theta(t)=\bar{P}(t)x(t)$, where
$\bar{P}(t)$ satisfies
\begin{eqnarray}
0=\dot{\bar{P}}(t)+\bar{P}(t)A_0(t)+\bar{P}(t)D_0\bar{P}(t)+A_0'(t)\bar{P}(t)+\Big(\bar{P}(t)B_0+C_0'(t)\Big)K_1(t).\label{new1}
\end{eqnarray}
In addition, substituting $\Theta(t)=\bar{P}(t)x(t)$ into (\ref{n2}) yields
\begin{eqnarray}
0=C_0(t)+B_0'\bar{P}(t).\label{new2}
\end{eqnarray}
Thus, we can reformulate (\ref{new1}) as
\begin{eqnarray}
0=\dot{\bar{P}}(t)+\bar{P}(t)A_0(t)+\bar{P}(t)D_0\bar{P}(t)+A_0'(t)\bar{P}(t).\nonumber
\end{eqnarray}
Comparing this with (\ref{d3}), we find that $\bar{P}(t)=P_1(t)$. Thus, (\ref{d4}) follows from (\ref{new2}) and (\ref{d5}) follows from $\Theta(T)=0$ and $\Theta(T)=P_1(T)x(T)$.

Second, the case where the controller $u_1(t)$ is of open--loop form can be solved similarly to the closed--loop case. This completes the proof.

\section{Proof of Theorem~\ref{Th-1}}\label{apendix3}

\emph{Proof of Sufficiency.} Under the condition (\ref{d4}), it is sufficient to verify that $P_1(T)x(T)=0$,
given Theorem~\ref{lemir}. To do this, we first state a formula relating $P_1(t)x(t)$ to the control $u_1(t)$ in terms of its dynamics.
Similar to (\ref{df2}), the dynamics of $P_1(t)x(t)$ are given by
\begin{eqnarray}
\frac{d[P_1(t)x(t)]}{dt}=-A_0'(t)P_1(t)x(t)-C_0'(t)u_1(t).\nonumber
\end{eqnarray}
Solving this differential equation yields
\begin{eqnarray}
P_1(t)x(t)
&=&\int_{t}^TP_2(t,s)C_0'(s)u_1(s)ds+P_2(t,T)C,  \label{df6}
\end{eqnarray}
where $C=P_1(T)x(T)$.

Next, we aim to prove that $C=0$ under the controller $u_1(t)$ defined in (\ref{df5}). If (\ref{df4}) holds, then for any $x_0$, there exists a $\zeta$ such that
$P_1(t_0)x_0=G_1[t_0,T]\zeta, $
where $\zeta=G_1^{\dag}[t_0,T]P_1(t_0)x_0$.
We can now rewrite $u_1(t)$ in (\ref{df5}) as
$u_1(t)=C_0(t)P_2'(t_0,t)\zeta$. By substituting $u_1(t)$ into (\ref{df6}), we obtain
\begin{eqnarray}
P_1(t_0)x_0
=\Big[\int_{t_0}^TP_2(t_0,s)C_0'(s)C_0(s)P_2'(t_0,s)ds\Big]\zeta+P_2(t_0,T)C=P_1(t_0)x_0+P_2(t_0,T)C, \nonumber
\end{eqnarray}
Since $P_2(t_0,T)$ is invertible, we have $C=0$, implying that $P_1(T)x(T)=0$. This completes the proof of sufficiency based on Theorem~\ref{lemir}.

\emph{Proof of Necessity.}  If the control problem is solvable, it follows from Theorem~\ref{lemir} that there exists a $P_1(t)$ such that (\ref{d4}) holds.
We now prove that (\ref{df4}) does indeed hold. Otherwise, we would have that $Range \big[P_1(t_0)\big]\nsubseteq
Range \Big(G_1[t_0,T]\Big)$, meaning that a non--zero vector $\rho$ would exist such that
$\rho'P_1(t_0)\rho\neq0,\rho'G_1[t_0,T]\rho=0$.
Then, we would obtain
\begin{eqnarray}
0=\rho'G_1[t_0,T]\rho=\rho'\Big[\int_{t_0}^TP_2(t_0,s)C_0'(s)C_0(s)P_2'(t_0,s)ds\Big]\rho
=\int_{t_0}^T\|C_0(s)P_2'(t_0,s)\rho\|^2ds,\nonumber
\end{eqnarray}
implying that
$C_0(s)P_2'(t_0,s)\rho=0$.
Thus, we would have
$\rho'\int_{t_0}^TP_2(t_0,s)C_0'(s)u_1(s)ds=0$.
Let $x_0=\rho$. From $\Theta(t_0)=P_1(t_0)x_0$, we would then have $\Theta(t_0)=P_1(t_0)\rho$.
Combining this with $\Theta(t)=\int_{t}^TP_2(t,s)C_0'(s)u_1(s)ds$ gives
$$\rho'P_1(t_0)\rho=\rho'\Big[\int_{t}^TP_2(t,s)C_0'(s)u_1(s)ds\Big]\rho=0.$$
This is a contradiction, so (\ref{df4}) must hold, completing the proof.

\section{Proof of Theorem~\ref{thd3}}\label{apendix4}

Let
\begin{eqnarray}
y(t)=\mathcal{T}_1'(t)x(t)=\left[
                           \begin{array}{c}
                             y_1(t) \\
                             y_2(t) \\
                           \end{array}
                         \right]. \label{Cls-1}
\end{eqnarray}
Then, using (\ref{n3}) and the feedback controller $u_1(t)=K(t)x(t)$, we have
\begin{eqnarray}
\dot{y}(t)&=&\dot{\mathcal{T}}_1'(t)x(t)+\mathcal{T}_1'(t)\Big(A_0(t)x(t)+D_0\Theta(t)+B_0K(t)x(t)\Big)\nonumber\\
&=&\Big[\dot{\mathcal{T}}_1'(t)\mathcal{T}_1(t)+\mathcal{T}_1'(t)\Big(A_0(t)+D_0P_1(t)+B_0K(t)\Big)\mathcal{T}_1(t)\Big]\mathcal{T}_1'(t)x(t)\nonumber\\
&=&\Big(\left[
     \begin{array}{c}
       \tilde{T}_{1}(t)\\
       \tilde{T}_{2}(t)\\
     \end{array}
   \right]+\left[
                                 \begin{array}{c}
                                   \hat{A}_{1}(t) \\
                                   \hat{A}_{2}(t) \\
                                 \end{array}
                               \right]+\left[
     \begin{array}{c}
       B_{1}(t)\\
       B_{2}(t)\\
     \end{array}
   \right]\mathcal{T}_1'(t) K(t)\mathcal{T}_1(t)\Big)y(t)\nonumber\\
&=&\left[
     \begin{array}{c}
       \tilde{T}_{1}(t)+\hat{A}_{1}(t)+B_{1}(t)\mathcal{T}_1'(t)K(t)\mathcal{T}_1(t)\\
       \tilde{T}_{2}(t)+\hat{A}_{2}(t)+B_{2}(t)\mathcal{T}_1'(t)K(t)\mathcal{T}_1(t)\\
     \end{array}
   \right]y(t).\nonumber
\end{eqnarray}
By applying (\ref{df7}), we obtain
$\left[
  \begin{array}{c}
    \dot{y}_1(t) \\
    \dot{y}_2(t) \\
  \end{array}
\right]=\left[
          \begin{array}{cc}
            \frac{I}{t-T} & 0 \\
            * & * \\
          \end{array}
        \right]\left[
                 \begin{array}{c}
                   y_1(t) \\
                   y_2(t) \\
                 \end{array}
               \right]$.
This implies that
$\dot{y}_1(t)=\frac{I}{t-T}y_1(t)$. Then, solving this equation gives us
$y_1(t)=\frac{T-t}{T-t_0}y_1(t_0)$,
further implying that
\begin{eqnarray}
y_1(T)=0.\label{d23}
\end{eqnarray}
Since
\begin{eqnarray}
0=\mathcal{T}_1'(T)P_1(T)\mathcal{T}_1(T)\mathcal{T}_1'(T)x(T)=\left[
         \begin{array}{cc}
           \hat{P}(T) & 0 \\
           0 & 0 \\
         \end{array}
       \right]y(T)=\hat{P}(T)y_1(T),
\end{eqnarray}
we can combine this with the invertibility of $\hat{P}(t)$ to obtain $y_1(T)=0$. Since, we also have (\ref{Cls-1}) and $y(T)=\mathcal{T}_1'(T)x(T)=\left[
                           \begin{array}{c}
                             0 \\
                             y_2(T) \\
                           \end{array}
                         \right]$, this yields
\begin{eqnarray}
P_1(T)x(T)&=&P_1(T)\mathcal{T}_1(T) \left[
                           \begin{array}{c}
                             0 \\
                             y_2(T) \\
                           \end{array}
                         \right] =\mathcal{T}_1(T)\mathcal{T}'_1(T) P_1(T)\mathcal{T}_1(T) \left[
                           \begin{array}{c}
                             0 \\
                             y_2(T) \\
                           \end{array}
                         \right]\nonumber \\
                         &=&\mathcal{T}_1(T) \left[
         \begin{array}{cc}
           \hat{P}(T) & 0 \\
           0 & 0 \\
         \end{array}
       \right]     \left[
                           \begin{array}{c}
                             0 \\
                             y_2(T) \\
                           \end{array}
                         \right] =0. \nonumber
 \end{eqnarray}
This completes the proof.

\section{Proof of Theorem~\ref{infir}}\label{apendix5}

\emph{Proof of Sufficiency.} Since solutions $P$ and $P_1$ exist to (\ref{dinf4}) and (\ref{d6}) that satisfy $P+P_1\geq 0$, we will show that
the cost function is bounded below by $x_0'(P+P_1)x_0$. By taking the derivative of $x'(t)(P+P_1)x(t)$, we obtain
\begin{eqnarray}
\frac{d}{dt}\Big(x'(t)(P+P_1)x(t)\Big)
&=&\Big(Ax(t)+Bu(t)\Big)'(P+P_1)x(t)+x'(t)(P+P_1)\Big(Ax(t)+Bu(t)\Big)\nonumber\\
&=&-x'(t)Qx(t)+x'(t)(P+P_1)BR^{\dag}B'(P+P_1)x(t)+u'(t)B'(P+P_1)x(t)\nonumber\\
&&+x'(t)(P+P_1)Bu(t),\nonumber
\end{eqnarray}
where we have used (\ref{dinf4}) and (\ref{d6}) to derive the last equality.
By integrating this from $0$ to $T$, we can further obtain
\begin{eqnarray}
&&\int_{0}^T\Big(x'(t)Qx(t)+u'(t)Ru(t)\Big)dt\nonumber\\
&=&x'(0)(P+P_1)x(0)-x'(T)(P+P_1)x(T)+\int_0^T\Big(u'(t)Ru(t)+x'(t)(P+P_1)BR^{\dag}B'(P+P_1)x(t)\nonumber\\
&&+u'(t)B'(P+P_1)x(t)+x'(t)(P+P_1)Bu(t)\Big)dt\nonumber\\
&=&x'(0)(P+P_1)x(0)-x'(T)(P+P_1)x(T)+\int_0^T\Big[\Big(u(t)+R^{\dag}B'(P+P_1)x(t)\Big)'R\Big(u(t)\nonumber\\
&&+R^{\dag}B'(P+P_1)x(t)\Big)+u'(t)(I-RR^{\dag})B'(P+P_1)x(t)+x'(t)(P+P_1)B(I-R^{\dag}R)u(t)\Big]dt\nonumber\\
&=&x'(0)(P+P_1)x(0)-x'(T)(P+P_1)x(T)+\int_0^T\Big(u(t)+R^{\dag}B'(P+P_1)x(t)\Big)'R\Big(u(t)\nonumber\\
&&+R^{\dag}B'(P+P_1)x(t)\Big)dt,\nonumber
\end{eqnarray}
where we have used $(I-RR^{\dag})B'(P+P_1)=0$ to derive the last equality, obtained from (\ref{d10}).
Since $u(t)\in \mathcal{U}$, we thus have $\lim_{T\rightarrow\infty}x'(T)(P+P_1)x(T)=0$. This implies that
\begin{eqnarray}
J(x_0;u)&=&\lim_{T\rightarrow\infty}\int_{0}^T\Big(x'(t)Qx(t)+u'(t)Ru(t)\Big)dt\nonumber\\
&=&x'(0)(P+P_1)x(0)+\lim_{T\rightarrow\infty}\int_0^T\Big(u(t)+R^{\dag}B'(P+P_1)x(t)\Big)'R\Big(u(t)+R^{\dag}B'(P+P_1)x(t)\Big)dt.\label{df10}
\end{eqnarray}
Because $R\geq 0$, we obtain $J(x_0;u)\geq x'(0)(P+P_1)x(0)$.

Next, we show that the controller (\ref{df8}) is stabilizing. Substituting (\ref{df8}) into (\ref{d1}) yields
\begin{eqnarray}
\dot{x}(t)&=&Ax(t)-BR^{\dag}B'(P+P_1)x(t)+BG_0Kx(t)\nonumber\\
&=&(A_0+D_0P_1)x(t)+BG_0Kx(t)\nonumber\\
&=&(A_0+D_0P_1)x(t)+BT_0^{-1}\left[
                               \begin{array}{c}
                                 0 \\
                                 K \\
                               \end{array}
                             \right]x(t).\label{df9}
\end{eqnarray}
Since $\Upsilon_{T_0}$ is of full row rank, there exists a $K_1$ such that $\Upsilon_{T_0}K_1=K$. From (\ref{jnYC1}), we find that
\begin{eqnarray}
T_0(I-R^{\dag}R)K_1=\left[
                              \begin{array}{c}
                                 0  \\
                                 \Upsilon_{T_0}\\
                              \end{array}
                            \right]K_1=\left[
                               \begin{array}{c}
                                 0 \\
                                 K \\
                               \end{array}
                             \right].\label{df11}
\end{eqnarray}
By substituting the above equation into (\ref{df9}), we then have
\begin{eqnarray}
\dot{x}(t)&=&(A_0+D_0P_1)x(t)+BT_0^{-1}T_0(I-R^{\dag}R)K_1x(t)\nonumber\\
&=&(A_0+D_0P_1)x(t)+B(I-R^{\dag}R)T_0^{-1}T_0(I-R^{\dag}R)K_1x(t)\nonumber\\
&=&(A_0+D_0P_1)x(t)+B_0\Upsilon_{T_0}K_1x(t)\nonumber\\
&=&(A_0+D_0P_1)x(t)+B_0Kx(t).\nonumber
\end{eqnarray}
Since $K$ was chosen, such that $A_0+D_0P_1+B_0K$ is stable, the above system is stable, and hence the controller (\ref{df8}) is stabilizing.

Finally, we substitute the stabilizing controller (\ref{df8}) into the cost function (\ref{df10}) to verify that (\ref{df8}) is an optimal controller, as desired. In fact,
with this controller, (\ref{df10}) becomes
\begin{eqnarray}
J(x_0;u)
&=&x'(0)(P+P_1)x(0)+\lim_{T\rightarrow\infty}\int_0^T\Big([-R^{\dag}B'(P+P_1)+G_0K]x(t)+R^{\dag}B'(P+P_1)x(t)\Big)'\nonumber\\
&&\times R\Big([-R^{\dag}B'(P+P_1)+G_0K]x(t)+R^{\dag}B'(P+P_1)x(t)\Big)dt\nonumber\\
&=&x'(0)(P+P_1)x(0)+\lim_{T\rightarrow\infty}\int_0^Tx'(t)K'G_0'RG_0Kx(t)dt.\label{df12}
\end{eqnarray}
By again using (\ref{df11}), it follows that $G_0K=T_0^{-1}\left[
                               \begin{array}{c}
                                 0 \\
                                 K \\
                               \end{array}
                             \right]=(I-R^{\dag}R)K_1$. This implies that $RG_0K=R(I-R^{\dag}R)K_1=0$.
Thus, with the controller (\ref{df8}), the cost function (\ref{df12}) reduces to
\begin{eqnarray}
J(x_0;u)=x'(0)(P+P_1)x(0).\nonumber
\end{eqnarray}
This shows that (\ref{df8}) is an optimal controller, and the optimal cost is $J^*(x_0;u)=x_0'(P+P_1)x_0$.
\\

\emph{Proof of Necessity.} Here, we derive the three conditions given in the theorem. First, we discuss
the results for the finite--horizon optimization problem. Considering the asymptotic behavior of the solutions to
the Riccati differential equations enables us to obtain the first and second conditions. Then, by applying the maximum principle, we find that the
stabilizability condition is as stated by the third condition. The detailed proof is given below.

First, based on Theorem~\ref{lemir}, there exists a $P^T(t)$ in (\ref{dr1}) and a $P_1^T(t)$ in (\ref{d3}) with terminal values of $P^T(T)=0$ and $P_1^T(T)$ such that
(\ref{d4}) holds, and there also exists a $u_{1}(t)$ that achieves (\ref{d5}), where $x(t)$ obeys (\ref{d24})
with initial value $x(0)=x_0$. In this case, the optimal cost is given by
\begin{eqnarray}
J_T^*(x_0;u)=x_0'\hat{P}^T(0)x_0,\label{d21}
\end{eqnarray}
where $\hat{P}^T(t)=P^T(t)+P_1^T(t)$. Given that $Q\geq0$ and $R\geq 0$, we have $J_T(x_0;u)\geq0$. Accordingly, for $T_1\leq T_2$, we obtain
$J_{T_1}(x_0;u)\leq J_{T_2}(x_0;u)$. Together with (\ref{d21}) and the arbitrariness of $x_0$, we thus find that
\begin{eqnarray}
\hat{P}^{T_1}(0)\leq \hat{P}^{T_2}(0).\label{d22}
\end{eqnarray}
In addition, consider the cost function
\begin{eqnarray}
J_T^t(x_0;u)=\int_{t}^T[x'(t)Qx(t)+u'(t)Ru(t)]dt.\nonumber
\end{eqnarray}
By applying a similar argument to that for Theorem~\ref{lemir}, the optimal cost yielded by minimizing $J_T^t(x_0;u)$ subject to (\ref{d1}) is given by
\begin{eqnarray}
J_T^t(x_0;u)=x'(t)\hat{P}^{T}(t)x(t).\nonumber
\end{eqnarray}
Since, for $t_1\leq t_2$, we have that
\begin{eqnarray}
J_T^{t_1}(x_0;u)\geq J_T^{t_2}(x_0;u),\nonumber
\end{eqnarray}
this implies that
\begin{eqnarray}
\hat{P}^{T}(t_1)\geq \hat{P}^{T}(t_2).\label{d35}
\end{eqnarray}
Combining (\ref{d22}) and (\ref{d35}), we see that $\hat{P}^T(t)$ is non--decreasing with respect to $T$ and that $\hat{P}^T(t)$ is non--increasing with respect to $t$.

Next, we show the uniform boundedness of $\hat{P}^T(t)$. Since there exists an optimal and stabilizing controller, there also exists a positive constant $c$ such that
\begin{eqnarray}
J_T^t(x_0;u)&\leq& \int_0^\infty\Big(x'(t)Qx(t)+u'(t)Ru(t)\Big)dt\leq c\|x_0\|^2.\nonumber
\end{eqnarray}
Combining this with (\ref{d21}), it follows that
$\hat{P}^T(t)\leq cI.$
As all the system matrices are time-invariant, $\hat{P}^T(t)$ is also time--invariant, i.e.,
$\hat{P}^T(t)=\hat{P}^{T-t}(0).$
Recalling (\ref{d22}) and (\ref{d35}), this shows that the limit $\lim_{T\rightarrow\infty} \hat{P}^T(t)=\hat{P}$ exists.
Moreover, by letting $t\rightarrow\infty$ in $\hat{P}^T(t)=P^T(t)+P_1^T(t)$, we see that $\hat{P}$ satisfies
\begin{eqnarray}
0&=&A'\hat{P}+\hat{P}A+Q-\hat{P}BR^{\dag}B'\hat{P}.\nonumber
\end{eqnarray}
This is exactly the same equation for $P$; hence, (\ref{dinf4}) is solvable. This further implies that (\ref{d6}) admits a solution $P_1$
and that $\hat{P}=P+P_1$.
Likewise, letting $t\rightarrow\infty$ in (\ref{d4}) yields
$C_0+B_0'P_1=0,$
which is exactly (\ref{d10}).

Finally, by applying the maximum principle, the optimal solution satisfies
\begin{eqnarray}
\dot{x}(t)&=&Ax(t)+Bu(t),\label{d27}\\
\dot{p}(t)&=&-A'p(t)-Qx(t),\label{d26}\\
0&=&Ru(t)+B'p(t),\label{d25}
\end{eqnarray}
with $\lim_{t\rightarrow\infty}p(t)=0$ and $x(0)=x_0$. Recalling that the optimal solution is also stabilizing, we obtain
\begin{eqnarray}
\lim_{t\rightarrow\infty}x(t)=0,\nonumber
\end{eqnarray}
and hence that
\begin{eqnarray}
\lim_{t\rightarrow\infty}Px(t)=0,\label{d29}\\
\lim_{t\rightarrow\infty}P_1x(t)=0.\label{d30}
\end{eqnarray}
Let
\begin{eqnarray}
p(t)=Px(t)+\Theta(t),\label{d23}
\end{eqnarray}
where $P(t)$ obeys equation (\ref{dinf4}) and $\Theta(t)$ is to be determined. From (\ref{d29}) and $\lim_{t\rightarrow\infty}p(t)=0$, we then have $\lim_{t\rightarrow\infty}\Theta(t)=0$.

By substituting (\ref{d23}) into (\ref{d25}), we obtain
\begin{eqnarray}
0&=&Ru(t)+B'Px(t)+B'\Theta(t).\nonumber
\end{eqnarray}
This implies that
\begin{eqnarray}
u(t)&=&-R^{\dag}\Big(B'Px(t)+B'\Theta(t)\Big)+(I-R^{\dag}R)z(t),\nonumber\\\label{d28}
\end{eqnarray}
and
\begin{eqnarray}
C_0x(t)+B_0'\Theta(t)=0,\label{d31}
\end{eqnarray}
where $z(t)$ is an arbitrary vector of compatible dimension.

Substituting (\ref{d28}) into (\ref{d1}) then reduces the state dynamics to
\begin{eqnarray}
\dot{x}(t)&=&Ax(t)-BR^{\dag}\Big(B'Px(t)+B'\Theta(t)\Big)+B(I-R^{\dag}R)z(t)\nonumber\\
&=&A_0x(t)+D_0\Theta(t)+B_0u_1(t),\label{d37}
\end{eqnarray}
where we have used $T_0(I-R^{\dag}R)z(t)=\left[
                                                                                              \begin{array}{c}
                                                                                                0 \\
                                                                                                \Upsilon_{T_0} \\
                                                                                              \end{array}
                                                                                            \right]z(t)\triangleq\left[
                                                                                              \begin{array}{c}
                                                                                                0 \\
                                                                                                u_1(t) \\
                                                                                              \end{array}
                                                                                            \right]
$ in the derivation of the last equality. Taking the derivative of (\ref{d23}) yields
\begin{eqnarray}
\dot{p}(t)=P\dot{x}(t)+\dot{\Theta}(t)=P\Big(A_0x(t)+D_0\Theta(t)+B_0u_1(t)\Big)+\dot{\Theta}(t).\nonumber
\end{eqnarray}
Comparing this with (\ref{d26}) and using (\ref{d6}), we have
\begin{eqnarray}
\dot{\Theta}(t)&=&-A_0'\Theta(t)-C_0'u_1(t).\label{d36}
\end{eqnarray}
We now prove that the solution to the FBDEs (\ref{d37}), (\ref{d36}), and (\ref{d31}) is
$\Theta(t)=P_1x(t)$, where $x(t)$ satisfies (\ref{d8}).  By taking the derivative of $P_1x(t)$, we obtain
\begin{eqnarray}
\frac{d}{dt}\Big(P_1x(t)\Big)
&=&P_1(A_0+D_0P_1)x(t)+P_1B_0u_1(t)\nonumber\\
&=&-A_0'P_1x(t)+P_1B_0u_1(t)\nonumber\\
&=&-A_0'P_1x(t)-C_0'u_1(t),\label{d38}
\end{eqnarray}
where we have used (\ref{d10}) to derive the last equality.
By comparing (\ref{d8}), (\ref{d38}), and (\ref{d10}) with (\ref{d37}), (\ref{d36}), and (\ref{d31}), we immediately obtain
\begin{eqnarray}
\Theta(t)=P_1x(t).\nonumber
\end{eqnarray}
Accordingly, the state dynamics are given by (\ref{d8}). To ensure the state's stability, these dynamics must be stabilizable, giving us the
third condition. This completes the proof.

\end{appendix}


\begin{thebibliography}{99}

\bibitem{Ho} Ho Y. Linear stochastic singular control problems. Journal of Optimization Theory and Application, 1972, 9(1): 24--31.

\bibitem{Speyer} Speyer J, Jacobson D. Necessary and sufficient conditions for optimality for singular control problems: a limit approach. Journal of Mathematical Analysis and Applications, 1971, 34(2): 239--266.

\bibitem{chenhanfu} Chen H-F. Unified controls applicable to general case under quadratic index. Acta Mathematicae Applicatae Sinica, 1982, 5(1): 45--52.

\bibitem{chenlizhou} Chen S, Li X, Zhou X. Stochastic linear quadratic regulators with indefinite control weight costs. SIAM Journal on Control and Optimization, 1998, 36(5): 1685--1702.

\bibitem{sunliyong} Sun J, Li X, Yong J. Open-loop and closed-loop solvabilities for stochastic linear quadratic optimal control problems. SIAM Journal on Control and Optimization, 2016, 54(5): 2274--2308.

\bibitem{Bell} Bell D J. Singular problems in optimal control-a survey. International Journal of Control, 1975, 21(2): 319--331.

\bibitem{clementsanderson} Clements D, Anderson B D O. Singular optimal control: The linear-quadratic problem. Springer-Verlag, New York, 1978.


\bibitem{Bellman} Bellman R, Glicksberg I, Gross O. Some aspects of the mathematical theory of control processes. Rand Corporation, R-313, 1958.

\bibitem{anderson} Anderson B D O, Moore J B. Optimal control: linear quadratic methods.
Englewood Cliffs, NJ: Prentice Hall, 1990.

\bibitem{lewis} Lewis F L, Vrabie D L, Syrmos V L. Optimal control.
John Wiley \& Sons, Inc., 2012.


\bibitem{SCIS1} Xu J, Shi J, Zhang H. A leader-follower stochastic linear quadratic differential game with time delay. SCIENCE CHINA Information Sciences, 2018, 61:112202.


\bibitem{SCIS2}	Qi Q, Zhang H. Time-inconsistent stochastic linear quadratic control for discrete-time systems. SCIENCE CHINA Information Sciences, 2017, 60:120204.

\bibitem{SCIS3}	Shi J, Wang G, Xiong J. Linear-quadratic stochastic Stackelberg differential game with asymmetric information. SCIENCE CHINA Information Sciences, 2017, 60:092202.

\bibitem{SCIS4} Ju P, Zhang H. Achievable delay margin using LTI control for plants with unstable complex poles. Science China Information Sciences, 2018, 61(9), 092203:1¨C092203:8.

\bibitem{hszhang} Zhang H, Li L, Xu J, Fu M. Linear quadratic regulation and stabilization of
discrete-time Systems with delay and multiplicative noise. IEEE Transactions on Automatic Control, 2015, 60(10): 2599--2613.


\bibitem{zhangIR} Zhang H, Xu J. On Irregular Linear Quadratic Control: Deterministic Case, arXiv:1711.09213.

\bibitem{zhangOpConsensus} Xu J, Zhang H. Consensus Control of Multi-agent Systems with Optimal Performance, arXiv:1803.09412.

\bibitem{hszhang1} Zhang H, Xu J. Control for It\^{o} stochastic systems with input delay. IEEE Transactions on Automatic Control, 2017, 62(1): 350--365.


\bibitem{Kalman} Kalman R E. Contributions to the theory of optimal control, Bol. Soc., Mat. Mexicana, 1960, 5: 102--119.


\bibitem{Letov} Letov A M. The analytical design of control systems. Automation \& Remote Control, 1961, 22: 363--372.


\bibitem{mp} Pontryagin L S, Boltyanskii V G, Gamkrelidze R V, Mishchenko E F. The mathematical theory of optimal process.
English translation. Interscience, 1962.


\bibitem{bellman} Bellman R. The theory of dynamic programming.
Bulletin of the American Mathematical Society, 60(6): 503--516, 1954.


\bibitem{Gurman} Gurman V. The method of multiple maxima and optimization problems for space maneuvers. Proc. Second Readings of K. E. Tsiolkovskii, Moscow, 1968, 39--51.

\bibitem{Moore} Moore J. The singular solutions to a singular quadratic minimization problem. International Journal of Control, 1974, 20(3): 383--393.


\bibitem{Williems} Willems J C, Kitapci A, Silverman L M. Singular optimal
control: a geometric approach. IAM Journal of Control and Optimization, 1986, 24(2): 323--337.


\bibitem{Gabasov} Gabasov R, Kirillova F M. High order necessary conditions for optimality. SIAM Journal on Control, 1972, 10: 127--168.


\bibitem{Krener}  Krener A J. The high order maximal principle and its application to singular extremals. SIAM Journal of Control and Optimization, 1977, 15: 256--293.


\bibitem{Hoehener} Hoehener D. Variational approach to second-order optimality conditions for control
problems with pure state constraints. SIAM Journal of Control and Optimization, 2012, 50: 1139--1173.


\bibitem{Bonnans} Bonnans J F, Silva F J. First and second order necessary conditions for
stochastic optimal control problems. Applied Mathematics \& Optimization, 2012, 65: 403--439.


\bibitem{XuZhang} Zhang H, Zhang X. Pointwise second-order necessary conditions for stochastic optimal controls, Part I: The case of convex control constraint.
SIAM Journal of Control and Optimization, 2015, 53(4): 2267--2296.


\bibitem{pinv} Penrose R. A generalized inverse of matrices, Mathematical Proceedings of the Cambridge Philosophical Society, 1955, 52: 17--19.

\bibitem{zhou2} Rami M Ait, Moore J B, Zhou X. Indefinite stochastic linear quadratic control
and generalized differential Riccati equation. SIAM Journal of Control and Optimization, 2002, 40: 1296--1311.



\bibitem{Kliger} Kliger I. Discussion on the stability of the singular trajectory with respect to ``Bang-Bang" control. IEEE Transactions on Automatic Control, 1964, 9(4): 583--585.

\bibitem{Hsia} Hsia T. On the existence and synthesis of optimal singular control with
quadratic performance index. IEEE Transactions on Automatic Control, 1967, 12(6): 778--779.

\bibitem{yong} Yong J,  Zhou X. Stochatic controls: Hamiltonian systems and HJB equations. Springer Verlag, 1999.

\bibitem{ZhangFangfang} Zhang F, Zhang H, Tan C, Wang W, et al. A new approach to distributed control for multi-agent systems based on approximate upper and lower bounds.
 International Journal of Control, Automation and Systems, 2017, 15(6): 2507--2515.


\bibitem{game} Xu J, Zhang H, Chai T. Necessary and sufficient condition for two-player Stackelberg strategy.
IEEE Transactions on Automatic Control, 2015, 60(5): 1356--1361.


\bibitem{meanfield} Zhang H, Qi Q, Fu M. Optimal stabilization control for discrete-time mean-field stochastic systems.
IEEE Transactions on Automatic Control, DOI: 10.1109/tac.2018.2813006.

\bibitem{ncs} Tan C, Zhang H. Necessary and sufficient stabilizing conditions for networked control
systems with simultaneous transmission delay and packet dropout.
IEEE Transactions on Automatic Control, 2017, 62(8): 4011--4016.


\end{thebibliography}
\end{document}